# Eratosthenes and Pliny, Greek geometry and Roman follies

Khristo N. Boyadzhiev
Department of Mathematics and Statistics
Ohio Northern University
Ada, OH 45810

**Abstract.** In this note we point out that supportive attitudes can bring to a blossoming science, while neglect and different values can quickly remove science from everyday life and provide a very primitive view of the world. We compare one important Greek achievement, the computation of the Earth meridian by Eratosthenes, to its later interpretation by the Roman historian of science Pliny.

## 1. Timeline

Something extraordinary happened in ancient Greece during the period 600-100 BC. During these five centuries the Greeks created the foundations of present-day science and mathematics. They achieved unprecedented results. The quality of these results is compatible only to those obtained seventeen centuries later. After the rise of the Roman Empire, however, the development of science and mathematics slowed down and came to an end. For several centuries after that there were practically no scientific activities in Europe.

Here is a short timeline.

**Pythagoras of Samos** (575-496 BC). Created the Pythagorean mathematical school. Obtained many theorems in Geometry and Number Theory.

**Socrates** (470–399 BC), a classical Greek philosopher, developed deduction, logic, and logical reasoning. His student Plato strongly influenced Aristotle, whom he taught on his part.

**Plato** (424-348 BC);

**Aristotle** (384-322 BC)

Aristotle worked with Alexander the Great and taught him science, mathematics and history.

**Alexander the Great** (356-323 BC) introduces Greek science and culture to his vast empire.

**Ptolemy I Soter** (367-283 BC), Alexander's general, ruler of Egypt, 323-283 BC. He started to build the most important scientific center of the ancient world - the Library of Alexandria (Bibliotheca Alexandrina).

**Euclid** (circa 300 BC); Author of the Elements.

**Ptolemy II Philadelphus** (309-246 BC), king of Egypt, 283-246 BC. Together with Ptolemy I he built the Musaeum at Alexandria, which included the famous Library of Alexandria.

**Eratosthenes of Cyrênê** (276-194 BC). Computed the length of the earth meridian and the circumference of the earth. Also known for the Sieve of Eratosthene.

**Archimedes of Syracuse** (287-212 BC).

**Heron of Alexandria** (about 10-70 AD). Engineer and mathematician.

**Pliny the Elder (**23-79 AD). A Roman historian and naturalist. Pliny is listed here for a special reason. His writings confirm the end of Greek science in the Roman empire. Scientific life proceeded to some extent in the Musaeum at Alexandria until the death of Hypatia.

**Hypatia of Alexandria** (350 – 415 AD). Possibly, the last prominent Greek scholar.

## 2. The educated rulers

The story of Eratosthenes is a logical part of the Greek history, a natural result of the general paradigm of this society – trust in science and in the scientific method. His work was partially made possible by the vision of Alexander the Great, who not only conquered new lands, but also promoted Greek culture and science. Alexander turned the small Egyptian port of Rhacotis into a new city, Alexandria, in order to make it the capital of his vast empire. Alexandria was well placed for trade with the Middle East. There was a little island in front of it, Pharos, to protect the city from the sea.

One of Alexander's generals, Ptolemy Soter, became king of Egypt (Ptolemy I) and put the foundations of the famous library. He also planned to build the seventh wander of the world – the lighthouse at Pharos (it was built in about 280 BC under his successor Ptolemy II). The lighthouse was a three-tier stone tower, approximately 120 meters high. On its top there was a platform where fires burn at night. The fires were reflected out to sea by metal mirrors. The lighthouse was named Pharos itself and today in many languages 'pharos' means a lighthouse.

Ptolemy Soter sponsored also Euclid in writing the Elements.

## 3. Eratosthenes and the meridian

Here we give a very brief account of how Eratosthenes computed the Earth meridian (and also the Earth circumference, radius, and diameter).

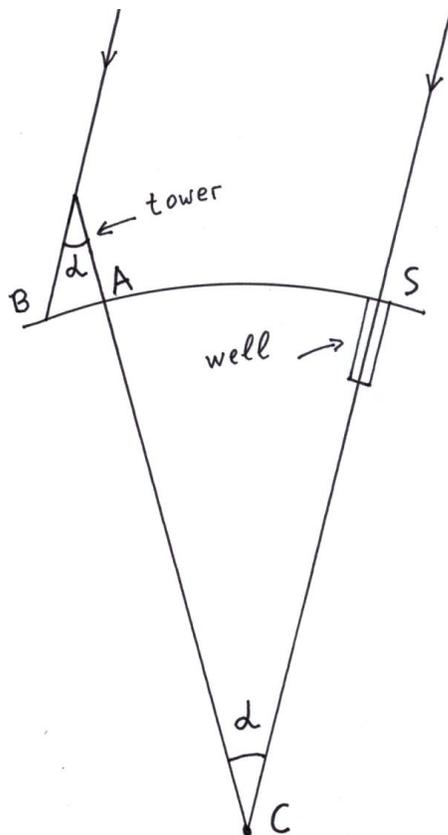

At the city of Syênè (Siena) in Egypt (today's name is Aswan), there was a deep vertical well, situated on the island of Elephantine in the Nile. People had noticed that every year on June 22 the sun could be seen from the bottom of the well. This meant the sun was exactly vertical. Also, at that day high objects had no shadow. Eratosthenes realized that on June 22 the sun, Siena and the center of the Earth were on one and the same line. He also realized that with some additional information this fact could be used to calculate the meridian of the Earth (and from this also the circumference of the Earth).

Convenient additional data came from Pharos in Alexandria. Eratosthenes measured the shadow AB cast by Pharos on June 22 and using the height of the tower (120 m) computed angle $\alpha$. It turned out to be 7 degrees and 12 minutes, exactly 1/50 of the complete angle 360 degrees. Eratosthenes also knew that when two parallel lines are cut by a third one, the corresponding angles at the intercepts are equal. He realized that angle $\alpha$ is the same as the angle at C between the two earth radii CA and CS reaching to Alexandria (A) and Siena (S) correspondingly.

Using proportions, Eratosthenes concluded that the distance AS between Alexandria and Siena is 1/50 of the complete circumference of the earth, the big circle passing trough these two cities. The distance AS was known to be about 5,000 stadia (most likely Egyptian stadia), which is about 787.5 km. Thus Eratosthenes computed the circumference to be 39,375 km, which is a very good result. (For a stydy on the exact value of $\alpha$ see [4]. The measurement of the meridian and Eratosthenes' method are discussed, for instance, in [1], [7], [10].)

It is good to mention here that the next in time reasonable measurement of the meridian was made by Caliph Al-Mamun in 820 AD, about ten centuries after Eratosthenes. Seven centuries later in 1525 came Fernal's result. Finally, DeLambre and Méchain in 1799 obtained an almost exact value. Their result was used by the French Academy of Sciences to give the meridional definition of the meter as one forty millionth part of the Earth circumference. Thus their value of the circumference was 40,000 km. After that, the definition of meter and the measurement of the meridian were slightly improved - the presently accepted value is slightly bigger [1].

In his fundamental study [5] Sir Thomas Heath writes: "*In truth all nations, in the West at all events, have been to school to the Greeks, in art, literature, philosophy, and science, the things which are essential to the rational use and enjoyment of human powers and activities, the things which make life worth living to a rational human being.*"([5, p.1]). Another researcher of the old Greeks, S. H. Butcher, explains their mindset: "*First, then, the Greeks, before any other people of antiquity, possessed the love of knowledge for its own sake. To see things as they really are, to discern their meanings and adjust their relations, was with them an instinct and a passion.*"([2, p.1]).

How did the old Greeks compare to the Romans? The mindset of the Romans was very different. They did not have a passion for science and love of knowledge – their passion was for military victories and governance. Philosophers (including mathematicians) and architects were often role models for the young Greeks. Socrates, for example, was a role model for many, including Plato. A legendary Greek hero and a role model, Odysseus, was not just a warrior – he was also an engineer and an inventor.

The situation in Rome was very different. Standard role models for the Roman youth were warriors, senators and patricians. The Romans were absorbed in a military culture and were busy conquering vast territories in Europe, Africa and Asia Minor. Their minds were tuned to the needs of the wars and the governance of the new territories. They needed mostly centurions and administrators.

**4. Pliny's version of the event**

Pliny the Elder (Plinii Secundi) was a Roman historian and naturalist who wrote a multivolume book Naturalis Historia, describing the world as the Romans knew it. He is the only known Roman writer who wrote about science. Pliny wrote about Eratosthenes too. We have included here some excerpts from his book [9] for illustration. The interesting phenomenon in Siena with the vertical well and objects with no shadow was known to Pliny, as seen from the following text

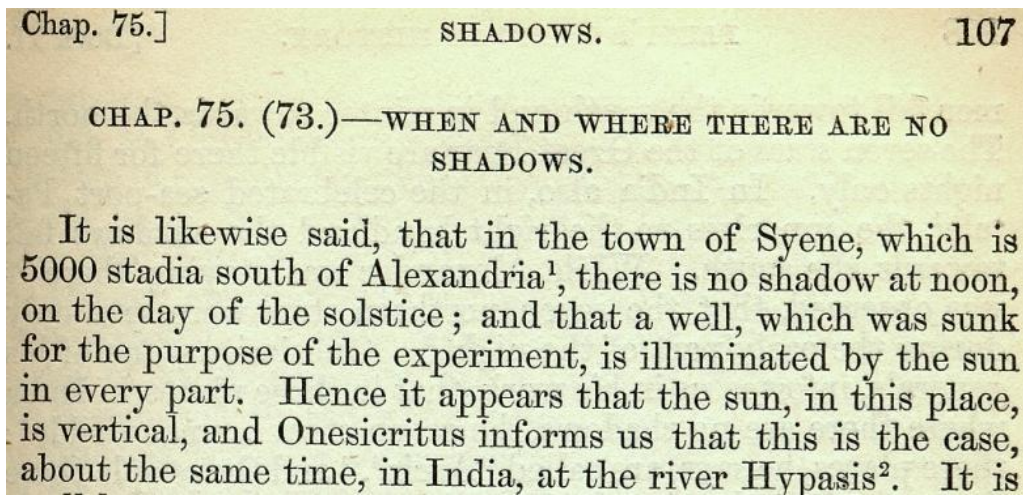

This observation, however, appears in a different part of the book and Pliny did not connect it to the work of Eratosthenes.

Let us see now how Pliny tells the story of Eratosthenes' computation of the meridian. Pliny writes, in fact, about the circumference of the Earth. This story can be found in the insert below, starting from "But Eratosthenes…" The reader will look in vain for some mathematical arguments and formulas.

merable islands lying off the coast of Germany[1], which have been only lately discovered.

The above is all that I consider worth relating about the length and the breadth of the earth[2]. But Eratosthenes[3], a man who was peculiarly well skilled in all the more subtle parts of learning, and in this above everything else, and a person whom I perceive to be approved by every one, has stated the whole of this circuit to be 252,000 stadia, which, according to the Roman estimate, makes 31,500 miles. The attempt is presumptuous, but it is supported by such subtle arguments that we cannot refuse our assent. Hipparchus[4], whom we must admire, both for the ability with which he controverts Eratosthenes, as well as for his diligence in everything else, has added to the above number not much less than 25,000 stadia.

(109.) Dionysodorus is certainly less worthy of confidence[5]; but I cannot omit this most remarkable instance of Grecian vanity. He was a native of Melos, and was celebrated for his knowledge of geometry; he died of old age in his native country. His female relations, who inherited his property, attended his funeral, and when they had for several successive days performed the usual rites, they are said to have found in his tomb an epistle written in his own name to those left above; it stated that he had descended from his tomb to the lowest part of the earth, and that it was a distance of 42,000 stadia. There were not wanting certain geometricians, who interpreted this epistle as if it had been sent from the middle of the globe, the point which is at the greatest distance from the surface, and which must necessarily be the centre of the sphere. Hence the estimate has been made that it is 252,000 stadia in circumference.

The remark about Dionysidorus is very interesting. Obviously, Eratosthenes' method, his "subtle parts of learning" and "subtle arguments" were beyond Pliny's understanding. Although he called him "a person, whom I perceive to be approved by everyone", Pliny himself did not seem to be convinced that just mental work and some simple observations could provide the length of the meridian and the radius of the Earth ("presumptuous attempt"). He felt the necessity to mention a fictitious voyage to the "middle of the globe" for directly measuring the radius. Anything more abstract seemed quite suspicious to the Roman mind that believed (and could understand) only direct measurement by hand. Pliny's attitude is very strong evidence that science (and geometry in particular) had no place in the powerful Roman Empire.

Writers on the history of mathematics usually give praises to the Greeks, but rarely comment on state of affairs in Rome. Some of them, however, do comment. Joseph Hofmann [6, p. 50] shortly states: "*The Romans produced no original mathematical accomplishment.*" Hofmann explains this partly with their awkward notational system (the Roman numbers). Florian Cajori writes in his classical work [3, p. 89] "*Although the Romans excelled in the science of government and war, in philosophy, poetry, and art they were mere imitators. In mathematics they did not even rise to the desire for imitation. … A science of geometry with definitions, postulates, axioms, rigorous proofs, did not exist there*".

One could ask: How was it possible for Rome to build the Coliseum, the aqueducts, and so many other remarkable constructions with such disregard of science and mathematics? The answer is simple and alas, gruesome – all Roman architecture was done by slaves, mostly enslaved Greek engineers. They worked from dawn to dusk, lived poorly, and were looked upon by the free Roman citizens. Under such conditions the Greeks could not maintain a vibrant scientific community and their knowledge as well as their skills slowly dissipated. Science activities continued for some time in Alexandria, but with the death of Hypatia in 415 AD they came to an end. The legacy of the Roman Empire includes several "dark" centuries for European science ([8], [11]).

(From the internet source http://www.roman-empire.net we read:
"*The state's public works were largely completed and maintained by slaves. Also, the government's state bureaucracy depended very much on educated slaves to keep the administration of the empire running. Even key institutions like the state's mints or the distribution of the corn dole to poor Romans depended on slaves. Other educated slaves also kept the private industries going, by functioning as their accountants and clerks.*

*Vital services were provided by literate slaves who served as teachers, librarians, scribes, artists and entertainers - even doctors.*")

Nothing illustrates better the contrast between the mindsets of the Greeks and the Romans, than the legend about the death of Archimedes. When Syracuse was overtaken by the Romans in 212 BC, a Roman soldier approached Archimedes who was intensely involved in some computations. On the question who he was, Archimedes simply answered: "Do not disturb me". Then the soldier killed him, just in case. The engraving below is after a painting by the French painter Gustave Courtois (1853-1923).

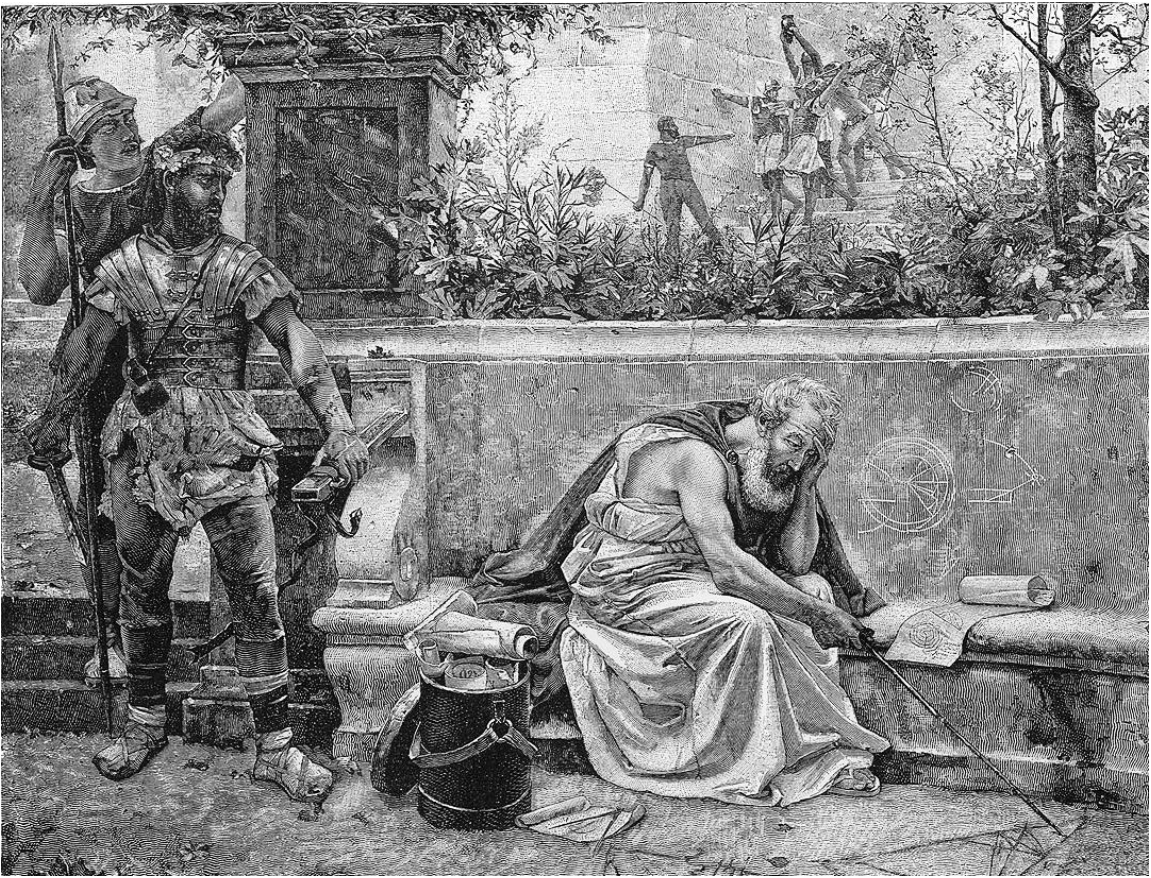